\documentclass{birkart}

 \newtheorem{thm}{Theorem}[section]
 \newtheorem{cor}[thm]{Corollary}
 \newtheorem{lem}[thm]{Lemma}
 \newtheorem{prop}[thm]{Proposition}
 \theoremstyle{definition}
 
 \theoremstyle{remark}
 \newtheorem{rem}[thm]{Remark}
 \newtheorem*{ex}{Example}
 \numberwithin{equation}{section}

\newcommand{\C}{{\mathbb C}}
\newcommand{\D}{{\mathbb D}}
\newcommand{\T}{{\mathbb T}}
\newcommand{\R}{{\mathbb R}}

\newcommand{\calH}{{\mathcal H}}
\newcommand{\calD}{{\mathcal D}}

\newcommand{\calE}{{\mathcal E}}
\newcommand{\calP}{{\mathcal P}}

\newcommand{\re}{\,{\rm Re}\,}

\newcommand{\diff}{{\mathrm d}}

\newcommand{\unitary}{{\,\mathcal U\,}}

\newcommand{\zed}{{\mathbf z}}
\newcommand{\wed}{{\mathbf w}}
\newcommand{\supp}{\,\mathrm{supp}\,}
\newcommand{\weight}{\Phi}

\newcommand{\corr}{\Gamma}
\newcommand{\Aset}{{\mathfrak A}}
\newcommand{\meas}{\sigma}
\newcommand{\Meas}{\Sigma}

\begin{document}
%
\title[Quantum Hele-Shaw flow]
 {
Quantum Hele-Shaw flow}
\author[Hedenmalm]{H\aa{}kan Hedenmalm}

\address{%
Department of Mathematics\\
Royal institute of Technology\\
S--100 44 Stockholm\\
Sweden}

\email{haakanh@math.kth.se}

\author[Makarov]{Nikolai Makarov}
\address{%
Department of Mathematics\\
California Institute of Technology\\
Pasadena, CA 91125\\
USA}
\email{makarov@caltech.edu}




\begin{abstract}
In this note, we discuss the quantum Hele-Shaw flow, a random measure 
process in the complex plane introduced by the physicists P.Wieg\-mann, 
A. Zabrodin, et al. This process arises in the theory of electronic 
droplets confined to a plane under a strong magnetic field, as well as 
in the theory of random normal matrices. We extend a result of Elbau 
and Felder \cite{ElFe} to general external field potentials, and  also show 
that if the potential is $C^2$-smooth, then the quantum Hele-Shaw flow 
converges, under appropriate scaling, to the  classical (weighted) 
Hele-Shaw flow, which can be modeled in terms of an obstacle problem.
\end{abstract}

\maketitle

\section{Introduction}

\noindent\bf A weighted distribution of $N$-tuples of complex numbers.
\rm
The intention of this paper is to provide a rigorous mathematical
treatment of the theoretical physics model due to Wiegmann, Zabrodin, et
al., which is concerned with the behavior of a finite number of electrons 
(charged fermions) in a strong magnetic field in the context of Quantum 
Theory. In part, this was done very recently by Elbau and Felder.
We will say more about this below.

Let $N$ be a positive integer, and introduce the notation
$$\zed_{[1,N]}=(z_1,z_2,z_3,\ldots,z_N)$$
for a point in $N$-dimensional complex Euclidean space $\C^N$. As a 
real vector space, we have the identification $\C^N\simeq\R^{2N}$. We let
$\lambda_{2N}$ be the standard Lebesgue measure on $\C^N$, so that in
particular,
$$\diff\lambda_{2N}\big(\zed_{[1,N]}\big)=\diff\lambda_{2}(z_1)\cdots
\diff\lambda_{2}(z_N).$$
The {\sl van der Monde determinant} is the quantity
$$\triangle\big(\zed_{[1,N]}\big)=\prod_{j,k:\,j<k}(z_j-z_k),$$ 
where we tacitly assume that the parameters $j$ and $k$ range over the 
set $\{1,\ldots,N\}$. Let $\weight:\C\to\R$ be a $C^2$-smooth function, 
which grows so quickly that 
\begin{equation}
\weight(z)\ge A\,\log|z|+O(1)
\label{est-below}
\end{equation}  
holds for all positive real constants $A$ (with $O(1)$ dependent on $A$).
Let $\beta$ denote a positive real parameter. We let 
$\mu_{\beta,\weight}$ be the finite positive Borel measure on $\C$
given by
\begin{equation}
\diff\mu_{\beta,\weight}(z)=
\exp\big\{-\beta\,\weight(z)\big\}\,\diff\lambda_2(z),
\qquad z\in\C,
\label{eq-mu}
\end{equation} 
We shall consider the probability measure
\begin{equation}
\diff\Pi_{N,\beta,\weight}\big(\zed_{[1,N]}\big)=
\frac1{Z_N}\,\big|\triangle\big(\zed_{1,N]}\big)\big|^2\,
\diff\mu_{\beta,\weight}(z_1)\cdots
\diff\mu_{\beta,\weight}(z_N),
\label{eq-distr-Pi}
\end{equation}
where $Z_N$ is the normalizing constant:
\begin{equation}
Z_N=\int_{\C^N}\big|\triangle\big(\zed_{1,N]}\big)\big|^2\,
\diff\mu_{\beta,\weight}(z_1)\cdots
\diff\mu_{\beta,\weight}(z_N).
\label{eq-distr-ZN}
\end{equation}
We intend to study the typical behavior of (\ref{eq-distr-Pi}) for large 
values of $N$. To make this more precise, let $n$ be a positive integer with
$1\le n\le N$, and consider the marginal distribution of 
$$\zed_{[1,n]}=(z_1,z_2,z_3,\ldots,z_n),$$
as $N\to+\infty$ while $n$ is kept fixed. This {\sl $n$-point marginal 
distribution} is given by the probability measure
\begin{multline}
\diff\Pi^{(n)}_{N,\beta,\weight}(\zed_{[1,n]})=
\frac1{Z_N}\,\Bigg\{\int_{\C^{N-n}}\big|\triangle(\zed_{[1,N]})\big|^2\,
\diff\mu_{\beta,\weight}(z_{n+1})\cdots
\diff\mu_{\beta,\weight}(z_{N})\Bigg\}\\
\times\diff\mu_{\beta,\weight}(z_1)\cdots\diff\mu_{\beta,\weight}(z_n).
\label{eq-distr-Pi2}
\end{multline}
The measure $\corr^{(n)}_{N,\beta,\weight}$ obtained by multiplying 
$\Pi^{(n)}_{N,\beta,\weight}$ by $N!/(N-n)!$ is called the {\sl $n$-point 
correlation measure}. The distribution (\ref{eq-distr-Pi})
appears from at least two different contexts: one is the theory of 
eigenvalues of random normal matrices, and the other a physical configuration 
of $N$ electrons localized to a two-dimensional plane, while being exposed to
a strong magnetic field (the strength is regulated by the
parameter $\beta$) perpendicular to the given plane. In the second instance, 
$\beta\weight$ is a scalar magnetic potential, $z_1,\ldots,z_N$ are the 
locations of the $N$ electrons, and one looks for time-independent solutions
to the Schr\"odinger equation. The resulting cloud of electrons is an
instance of the Aharonov-Bohm effect, where the magnetic potential rather
than the magnetic field manifests itself. For details, we ask the reader to 
consult the work of Wiegmann, Zabrodin, et al. \cite{ABWZ}, 
\cite{KKMWWZ}, \cite{MWZ}, \cite{KMZ}, and \cite{Wieg}. In the first 
instance (random normal matrices), we should mention first that normal 
$N\times N$ matrices have the decomposition $M=U^*DU$, where $D$ is 
diagonal (consisting of the eigenvalues), and $U$ is unitary, that is, 
$U^*U=UU^*=I$, where $I$ is the identity matrix. To focus on the eigenvalues, 
then, one should remove as much as possible of the unitary part; in other 
words, we should mod out with respect to 
$\unitary(N)/\unitary_{\mathrm d}(N)$, where $\unitary(N)$
is the group of unitary $N\times N$ matrices, and $\unitary_{\mathrm d}(N)$  
is the subgroup of unitary diagonal matrices. If this is done carefully, 
and a weight is introduced to get a probability measure on the normal 
matrices, we get a joint distribution of the eigenvalues of the type 
(\ref{eq-distr-Pi}). For details of this, see \cite{ElFe} and \cite{ChZa};
see also Mehta's book \cite{Meh}.  

The standard models for the eigenvalue distribution of self-adjoint matrices 
are by now quite well-known (see \cite{Meh}, \cite{Kurt}). In the limit, 
we get such models from our probability measure (\ref{eq-distr-Pi}) by 
letting $\weight(z)$ tend to $+\infty$ on $\C\setminus\R$. For instance, the
semicircle law of Wigner can be obtained by taking such a limit with the
appropriate interpretation.  

\medskip

\noindent\bf Quantum Hele-Shaw flow vs Hele-Shaw flow.
\rm We show that the $n$-point correlation measure 
$$\corr^{(n)}_{N,\beta,\weight}=\frac{N!}{(N-n)!}\,
\Pi^{(n)}_{N,\beta,\weight}$$
is given in terms of a determinantal process involving reproducing kernels
of polynomial subspaces of a Bargmann-Fock type space, and that it defines 
a monotonic growth process in $N$ (for fixed $\beta$ and $\weight$): 
$$\corr^{(n)}_{N,\beta,\weight}\le\corr^{(n)}_{N+1,\beta,\weight},\qquad
N=1,2,3,\ldots.$$
For this reason, and for reasons which will be made clearer below, we
call this {\sl Quantum Hele-Shaw flow}. The analogy with classical
Hele-Shaw flow was made earlier by Wiegmann and Zabrodin. The ratio
$N/\beta$ seems like a natural growth parameter for the process. For
mathematical reasons, however, it is more natural to use
\begin{equation}
\tau=\frac{N-1}{\beta}
\label{def-tau}
\end{equation}
instead. 

The weighted energy of a compactly supported Borel probability measure 
$\meas$ is (for a given positive parameter $\tau$) 
\begin{equation*}
\calE_{\weight}[\meas]\left(\frac{1}{\tau}\right)
=\frac1{2\tau}\iint_{\C^2}\bigg\{2\tau\log\frac{1}{|z-w|}+\weight(z)+
\weight(w)\bigg\}
\,\diff\meas(z)\diff\meas(w);
\end{equation*} 
the unique compactly supported Borel probability measure which minimizes
this energy is called the {\sl extremal measure}, or 
{\sl equilibrium measure}, and denoted by $\widehat\meas_{\weight,\tau}$, 
or simply by $\widehat\meas$, if no misunderstanding is possible. 
The existence and uniqueness of the extremal measure are treated in Chapter
1 of Saff's and Totik's book \cite{SaTo}. 

Our first main result reads as follows.

\begin{thm}
If $\tau$, $0<\tau<+\infty$, is kept fixed, and $\beta,\,N$ are connected via
(\ref{def-tau}), then, as $N\to+\infty$, the marginal distribution has the 
weak-star limit
$$\diff\Pi^{(n)}_{N,\beta,\weight}(\zed_{[1,n]})\to
\diff\widehat\meas(z_1)\cdots\diff\widehat\meas(z_n).$$  
\label{thm-I}
\end{thm}

This theorem is essentially as in Johansson's paper \cite{Kurt}, but with 
important additional technical details. Elbau and Felder \cite{ElFe}
were the first to realize that Johansson's argument extends to the case of 
normal random matrices. The present version is more general -- Elbau and 
Felder required that the potential be infinite in a neighborhood of infinity.

We discuss the extremal measure in Section \ref{sec2}, and connect it to an
obstacle problem in Section \ref{sec3}. This leads to our second main theorem.
For the formulation, we need the following standard notation: $\Delta$ is
the standard Laplacian in the plane, while given a subset $E$ of the complex 
plane $\C$, $1_E$ is the characteristic function of the set $E$.

\begin{thm}
For each $\tau$, $0<\tau<+\infty$, there exists an unbounded open set 
$\calD_\weight(\tau)$ in $\C$ whose complement is compact, such that the 
extremal measure has the form
$$\diff\widehat\meas_{\weight,\tau}(z)=\frac1{4\pi\tau}\,
1_{\C\setminus\calD_\weight(\tau)}(z)\,\Delta\weight(z)\,\diff\lambda_2(z),
\qquad z\in\C.$$
Moreover,  as $\tau$ increases, the open sets $\calD_\weight(\tau)$ form a
continuous decreasing chain, and the measures 
$\tau\,\widehat\meas_{\weight,\tau}$ grow with $\tau$.
\label{thm-II}
\end{thm}

This theorem justifies the claim of Wiegmann, Zabrodin, et al., that in the 
limit one gets a monotone family of domains (and not only measures) in the 
complex plane. The proof is an application of the general theory of 
free boundary problems due to Caffarelli et al.

We should add that, if the weight function $\weight$ is real analytic, then 
the boundary $\partial\calD_\weight(\tau)$ consists of real analytic curves, 
except for finitely many cusps (pointing outward) or so-called contact 
points. The evolution of $\calD_\weight(\tau)$ is quite similar to the 
Hele-Shaw flow on hyperbolic or weakly hyperbolic surfaces described in 
\cite{HS} and \cite{HO}. In fact, the present situation corresponds to 
Hele-Shaw flow on surfaces where the metric may be negative in portions of 
space, since $\Phi$ is a kind of metric potential, with associated metric  
$$\diff s(z)^2=\Delta\Phi(z)\,|\diff z|^2.$$
Generally speaking, Hele-Shaw flow models the expansion of viscous fluid in a
surrounding medium (which we may think of as vacuum), as new fluid is 
injected at a constant rate at one or several source points. In the present
setting, we have a single source point at infinity.

A final remark we wish to make here is that essentially all results of
the paper remain valid with little change if the $C^2$-smoothness of $\weight$
is weakened to local $W^{2,p}$-smoothness, for some $p$, $1<p<+\infty$. 

\section{Quantum Hele-Shaw flow and Bargmann-Fock type spaces}

\noindent\bf Weighted Bargmann-Fock spaces and reproducing kernels. \rm
Let $\mu$ be a positive finite Borel measure on the complex plane $\C$, with
$$0<\int_\C |z|^{2n}\diff\mu(z)<+\infty,\qquad n=0,1,2,3,\ldots.$$
The Hilbert space $L^2(\C,\mu)$ consists of the square-integrable 
complex-valued functions with respect to $\mu$, where two functions are 
identified if the coincide $\mu$-almost everywhere. We shall assume that $\mu$
is supported on an infinite set, so that $L^2(\C,\mu)$ becomes 
infinite-dimensional. The above integrability condition ensures that all 
(analytic) polynomials belong to $L^2(\C,\mu)$. 
Let $P^2(\C,\mu)$ denote the closure of the linear space of polynomials in 
$L^2(\C,\mu)$. Also, for $N=1,2,3,\ldots$, let $P^2_N(\C,\mu)$ be the 
$N$-dimensional linear space of polynomials of degree $\le N-1$, supplied 
with the Hilbert space structure of $L^2(\C,\mu)$. The point evaluations
in $\C$ are linear functionals on $P^2_N(\C,\mu)$, and thus given by the
inner product with an element of the space:
$$p(z)=\langle p,K_\mu(\cdot,z)\rangle_{L^2(\mu)}=\int_\C
\bar K_{\mu,N}(w,z)\,p(w)\,\diff\mu(w),\qquad z\in\C,$$   
where, for fixed $z\in\C$, the function $K_{\mu,N}(\cdot,z)$ is an element
of $P^2_N(\C,\mu)$, and so is $p$. The function $K_{\mu,N}$ is called the
{\sl reproducing kernel} for the space $P^2_N(\C,\mu)$. The monomials 
\begin{equation}
e_j(z)=z^{j-1},\qquad j=1,2,\ldots,N,
\label{mon-basis}
\end{equation}
form a vector space basis in $P^2_N(\C,\mu)$, but they need not be
orthogonal to one another. To rectify this, we may apply the Gram-Schmidt
process: let $\varphi_1=e_1$, and for $k=2,3,\ldots,N$, put
\begin{equation}
\varphi_k=e_k+\sum_{j=1}^{k-1}\lambda_{j,k}\,e_j,
\label{Gram-Schmidt}
\end{equation}
where the complex scalars $\lambda_{j,k}$ are chosen so that $\varphi_k$
becomes perpendicular to all the vectors $e_1,e_2,\ldots,e_{k-1}$. 
The functions $\varphi_1,\varphi_2,\varphi_3,\ldots$ obtained in this fashion
form an orthogonal sequence. To get an orthonormal sequence, we normalize:
$$\phi_j(z)=\frac{\varphi_j(z)}{\|\varphi_j\|_{L^2(\mu)}},\qquad z\in\C.$$
The reproducing kernel may then be conveniently expressed in terms of these 
orthogonal functions:
\begin{equation}
K_{\mu,N}(z,w)=\sum_{j=1}^N \phi_j(z)\,\bar\phi_j(w),\qquad 
(z,w)\in\C\times\C.
\label{kernelexp}
\end{equation}
\medskip

\noindent\bf Matrix-valued reproducing kernels and correlation measures. 
\rm For two vectors 
$$\zed_{[1,n]}=(z_1,\ldots,z_n)\quad\text{ and }\quad 
\wed_{[1,n]}=(w_1,\ldots,w_n),$$ 
we form the $n\times n$ matrix 
\begin{equation}
{\mathbf K}^n_{\mu,N}(\zed_{[1,n]},\wed_{[1,n]})=
\big[K_{\mu,N}(z_j,w_k)\big]_{j,k=1}^n,
\label{matrix-K}
\end{equation} 
which we may think of as a matrix-valued reproducing kernel. 

We return to the joint distribution function (\ref{eq-distr-Pi}), in the
slightly altered form 
\begin{equation}
\diff\Pi_{\mu,N}(\zed_{[1,N]})=
\frac1{Z_N}\,\big|\triangle(\zed_{[1,N]})\big|^2\,\diff\mu(z_1)\cdots
\diff\mu(z_N),
\label{distr-pi-1}
\end{equation}
where $\zed_{[1,N]}=(z_1,\ldots,z_N)$, and $Z_N$ is the normalizing constant: 
\begin{equation}
Z_N=
\int_{\C^N}\big|\triangle(\zed_{[1,N]})\big|^2\,\diff\mu(z_1)\cdots 
\diff\mu(z_N).
\label{normconst}
\end{equation}
We fix an integer $n$, $1\le n\le N$, and split 
$\zed_{[1,N]}=(\zed_{[1,n]},\zed_{]n,N]})$, with 
$\zed_{[1,n]}=(z_1,\ldots,z_n)$ and $\zed_{]n,N]}=(z_{n+1},\ldots,z_N)$. The 
{\sl $n$-point marginal distribution measure} is
\begin{multline}
\diff\Pi^{(n)}_{\mu,N}(\zed_{[1,n]})=\frac1{Z_N}\,
\Bigg\{\int_{\C^{N-n}}\big|\triangle(\zed_{[1,N]})\big|^2\,
\diff\mu(z_{n+1})\cdots\diff\mu(z_{N})\Bigg\}\\
\times\diff\mu(z_1)\cdots\diff\mu(z_n),
\label{distr-pi-2.1}
\end{multline}
while the {\sl $n$-point correlation measure} is 
\begin{equation}
\diff\corr^{(n)}_{\mu,N}\big(\zed_{[1,n]}\big)=\frac{N!}{(N-n)!}\,
\diff\Pi^{(n)}_{\mu,N}(\zed_{[1,n]}),\qquad \zed_{[1,n]}\in\C^n,
\label{distr-pi-2.2}
\end{equation}
where the integration takes place with respect to the variables 
$z_{n+1},\ldots,z_N$ only. For $n=1$, the measure $\corr^{(1)}_{\mu,N}$ 
describes the joint density of the eigenvalues. For measures $\mu$ of the
form (\ref{eq-mu}), with $\Delta\Phi$ a positive constant in a large disk 
centered at the origin, $\corr^{(1)}_{\mu,N}$ is rather similar to the 
majorization function of Aleman, Richter, and Sundberg \cite{ARS}; one major 
difference, though, is that here, we consider backward shift invariant 
subspaces ($P^2_N(\C,\mu)$), while in \cite{ARS}, the object of study are 
the forward shift invariant subspaces. 

\begin{prop}
We have that
$$\det\Big[{\mathbf K}^N_{\mu,N}(\zed_{[1,N]},\zed_{[1,N]})\Big]=
\frac{N!}{Z_N}\,\big|\triangle(\zed_{[1,N]})\big|^2.$$
\label{prop-n=N}
\end{prop}

\begin{proof}
We first apply simple row or column operations on the representation of the 
van der Monde determinant as the determinant of the $N\times N$ matrix with 
entries $e_j(z_k)$, for $j,k=1,\ldots,N$. The result is that the van der 
Monde determinant may be calculated based on the entries $\varphi_j(z_k)$ 
instead, whence the result follows from the formula defining the determinant 
in terms of permutations, as in \cite{Meh}, pp. 89--95.
\end{proof}

The following ``linear algebra'' description of $\corr^{(n)}_{\mu,N}$
in terms of the determinant of a matrix reproducing kernel is essentially
known (see \cite{Meh} for the Gaussian situation). 

\begin{thm}
We have that
$$\diff\corr^{(n)}_{\mu,N}(\zed_{[1,n]})=\det
\Big[{\mathbf K}^n_{\mu,N}(\zed_{[1,n]},\zed_{[1,n]})\Big]\,\,
\diff\mu(z_1)\cdots\diff\mu(z_n).$$
\label{thm-Lambda}
\end{thm}

\begin{proof}
We are to show that
$$\frac{N!}{(N-n)!}\,
\int_{\C^{N-n}}\big|\triangle(\zed_{[1,N]})\big|^2\,
\diff\mu(z_{n+1})\cdots\diff\mu(z_{N})
=Z_N\,\det\big[{\mathbf K}^n_{\mu,N}(\zed_{[1,n]},\zed_{[1,n]})\big].$$
We follow the approach of \cite{Meh}, pp. 89--95, and use that by Proposition
\ref{prop-n=N}, the identity is valid for $n=N$, which we take as our
starting point. Then successive applications of Theorem 5.2.1 in \cite{Meh} 
(suitable modified by including complex conjugates where needed) yields the 
result for general $n=1,\ldots, N$.
The proof is complete.
\end{proof}

\begin{cor}
The $n$-point correlation measure $\diff\corr^{(n)}_{\mu,N}$ increases
monotonically with increasing $N$. Its total mass equals $N!/(N-n)!$. 
\end{cor}

\begin{proof} That the total mass of $\diff\corr^{(n)}_{\mu,N}$ is
$N!/(N-n)!$ follows from the definition (\ref{normconst})
of the normalizing constant $Z_N$.
In view of (\ref{kernelexp}), we have
$$K_{\mu,N+1}(z,w)=K_{\mu,N}(z,w)+\phi_{N+1}(z)\,\bar\phi_{N+1}(w),$$
so that
$${\mathbf K}^n_{\mu,N+1}(\zed_{[1,n]},\wed_{[1,n]})=
{\mathbf K}^n_{\mu,N}(\zed_{[1,n]},\wed_{[1,n]})+
{\boldsymbol\Phi}_{\mu,N+1}(\zed_{[1,n]},\wed_{[1,n]}),$$
where 
$${\boldsymbol\Phi}_{\mu,N+1}(\zed_{[1,n]},\wed_{[1,n]})=
\big(\phi_{N+1}(z_j)\,\bar\phi_{N+1}(w_k)\big)_{j,k=1}^{n}.$$ 
The matrix ${\boldsymbol\Phi}_{\mu,N+1}(\zed_{[1,n]},\zed_{[1,n]})$ is 
positive definite, and so is 
${\mathbf K}^n_{\mu,N}(\zed_{[1,n]},\zed_{[1,n]})$. As the
determinant operation is monotonic with respect to increasing positive 
definiteness (see Proposition \ref{prop-A+B} below), 
the assertion follows from Theorem \ref{thm-Lambda}.
The proof is complete.
\end{proof}
\medskip

\noindent\bf A linear algebra result. 
\rm The following proposition is known; for instance, it should be possible 
to derive it from the Horn inequalities, which describe the possible 
eigenvalues of the sum of two Hermitian matrices in terms of the eigenvalues 
of the summands (see, e. g., the survey paper \cite{Ful}). However, we find 
it convenient to supply instead a short direct proof.

\begin{prop}
Let $A$ and $B$ be two positive semi-definite (Hermitian) $N\times N$ 
matrices. Then
$$\max\big\{\det(A),\det(B)\big\}\le\det(A+B).$$
\label{prop-A+B}
\end{prop} 

\begin{proof}
Clearly, by symmetry, it suffices to show that
$$\det(A)\le\det(A+B).$$
If $\det(A)=0$, the assertion follows trivially, as the determinant
of a positive semi-definite matrix is always $\ge0$. In the non-singular
case $\det(A)>0$, we note that $A^{-1/2}$ is a well-defined positive
semi-definite matrix, and that $A^{-1/2}\,B\,A^{-1/2}$ is positive 
semi-definite as well. We write
$$\det(A+B)=\det(A)\,\det\big(I+A^{-1/2}\,B\,A^{-1/2}\big),$$
where $I$ denotes the $N\times N$ identity matrix. The assertion now 
follows from the fact that the eigenvalues of $A^{-1/2}\,B\,A^{-1/2}$ are
all positive, so that the eigenvalues of $I+A^{-1/2}\,B\,A^{-1/2}$ are
all greater than $1$ (after all, the determinant equals the product of 
the eigenvalues). 
\end{proof}

\section{Weighted Fekete points and the extremal measure}\label{sec2}

\noindent\bf Weighted Fekete points. \rm
Let us analyze the probability measure (\ref{eq-distr-Pi}), written in the 
form
\begin{equation}
\diff\Pi_{N,\beta,\weight}(\zed_{[1,N]})=
\frac1{Z_N}\exp\bigg\{-\beta\sum_{j=1}^N \weight(z_j)\bigg\}\,
\big|\triangle(\zed_{[1,N]})\big|^2\,\diff \lambda_{2N}(\zed_{[1,N]}),
\label{distr-pi-3}
\end{equation}
where $\zed_{[1,N]}=(z_1,\ldots,z_N)$, and $Z_N$ is the normalization 
constant. The most likely configuration $\zed_{[1,N]}$ is the one that 
maximizes
$$\exp\bigg\{-\beta\sum_{j=1}^N \weight(z_j)\bigg\}\,
\big|\triangle(\zed_{[1,N]})\big|^2$$
among all vectors in $\C^N$. Considering that
$$\sum_{j,k:\,j<k}\big[\weight(z_j)+\weight(z_k)\big]=(N-1)\sum_{j=1}^N
\weight(z_j),$$
we may rewrite this expression:
\begin{multline}
\exp\bigg\{-\beta\sum_{j=1}^N \weight(z_j)\bigg\}\,
\big|\triangle(\zed_{[1,N]})\big|^2\\
=\prod_{j,k:\,j<k}
\Bigg\{\exp\bigg(-\frac{\beta}{N-1}\,\big[\weight(z_j)+\weight(z_k)\big]
\bigg)\,|z_j-z_k|^2\Bigg\}.
\label{eq-dens}
\end{multline}
For positive real $\theta$, let $E_\weight(z,w;\theta)$ denote the function
\begin{equation}
E_\weight(z,w;\theta)=|z-w|\,\exp\Big(-\theta\,\big[\weight(z)+
\weight(w)\big]\Big);
\label{E-fn}
\end{equation}
we are interested in maximizing
$$\prod_{j,k:\,j<k}E_\weight(z_j,z_k;\theta),$$
with 
$$\theta=\frac{\beta}{2(N-1)}.$$ 
Let $M_{\weight,N}(\theta)$ be this maximum raised to the power $2/(N(N-1))$:
$$M_{\weight,N}(\theta)=\sup_{\zed_{[1,N]}\in\C^N}\Bigg\{
\prod_{j,k:\,j<k}E_\weight(z_j,z_k;\theta)\Bigg\}^{2/(N(N-1))}.$$
The points $\zed_{[1,N]}=\zed_{[1,N]}^*$ which achieve this maximum are
called {\sl weighted Fekete points} (see \cite{SaTo}). In order that
the maximum be assumed, we should add an assumption on the weight $\weight$:
\begin{equation}
E(z,w;\theta)\to0\qquad \text{as}\quad\max\{|z|,|w|\}\to+\infty.
\label{reg-cond}
\end{equation}
An easy argument shows that (\ref{reg-cond}) follows from (\ref{est-below})
with an appropriate choice of the parameter $A$. 
It follows from Theorem 1.1 \cite[p. 143]{SaTo} that $M_{\weight,N}(\theta)$ 
{\sl is decreasing in} $N$, provided that $\weight$ is fixed and $\theta$ is 
held constant.
Moreover, under the same conditions, $M_{\weight,N}(\theta)\to 
M_{\weight}(\theta)$ as $N\to+\infty$, where 
\begin{equation}
M_{\weight}(\theta)=\exp(-\kappa_\weight(\theta)),
\label{limit-eq}
\end{equation}
and $\kappa_\weight(\theta)$ is results from a certain minimization problem,
outlined below. The quantity $M_\weight(\theta)$ can be thought of as a 
weighted capacity (of the whole plane).  
\medskip

\noindent\bf The extremal measure. \rm
Let $\calP_c(\C)$ denote the convex ``body'' of all compactly supported 
Borel probability measure $\meas$ on $\C$. For $\meas\in\calP_c(\C)$, we 
introduce the {\sl energy functional}
\begin{multline}
\calE_{\weight}[\meas](\theta)=2\theta\int_\C \weight(z)\,\diff\meas(z)+
\iint_{\C^2}\log\frac{1}{|z-w|}
\diff\meas(z)\diff\meas(w)\\
=\iint_{\C^2}\bigg\{\log\frac{1}{|z-w|}+\theta\,\weight(z)+
\theta\,\weight(w)\bigg\}
\,\diff\meas(z)\diff\meas(w)\\
=\iint_{\C^2}\log \frac{1}{E_\weight(z,w;\theta)}\,
\diff\meas(z)\diff\meas(w),
\label{energy}
\end{multline} 
and consider the problem of minimizing $\calE_{\weight}[\meas](\theta)$ over
all $\meas$. Define 
\begin{equation}
\kappa_\weight(\theta)=\inf_{\meas\in\calP_c(\C)}
\,\calE_{\weight}[\meas](\theta);
\label{cap-def}
\end{equation}
this quantity is related to $M_\weight(\theta)$ via (\ref{limit-eq}), 
in view of Theorem 1.3 \cite[p. 145]{SaTo}. Any $\meas\in\calP_c(\C)$
which achieves the maximum in (\ref{cap-def}) is called an {\sl extremal
(or equilibrium) measure}. In view of (\ref{reg-cond}), which follows from
(\ref{est-below}), it follows from Theorem 1.3 of \cite[p. 27]{SaTo} that 
there {\sl exists a unique extremal measure, which we will denote by } 
$\widehat\meas$, or, when
it necessary to indicate the dependence on the parameter $\theta$ as well
as on the weight $\weight$, by  $\widehat\meas_{\weight,\theta}$. 

\section{The extremal measure and Hele-Shaw flow}\label{sec3}

\noindent\bf A change of parameters. \rm
For reasons of convenience, we will now use $\tau=1/(2\theta)$ in place of 
$\theta$ as a positive real parameter.
As before, $\weight:\C\to\R$ is $C^2$-smooth, and subject to the growth 
requirement (\ref{est-below}) for every choice of the real parameter $A$. 
In particular, (\ref{reg-cond}) holds with $\theta=1/(2\tau)$.
As we mentioned in the previous section, by Theorem 1.3 \cite[p.27]{SaTo}, 
there exists a unique extremal measure $\widehat\meas\in\calP_c(\C)$ which 
achieves the infimum in (\ref{cap-def}). For $\widehat\meas$, the 
associated logarithmic potential
$$L\big[\widehat\meas\big](z)=
\int_\C\log\frac{1}{|z-w|}\,\diff\widehat\meas(w),
\qquad z\in\C,$$
is locally bounded in $\C$ (Theorem 4.3 \cite[p. 51]{SaTo}). So, for 
instance, the extremal measure $\widehat\meas$ has no point masses. 
To get an understanding of this extremal measure, we make a detour to
obstacle problems and Hele-Shaw flow. A general reference for Hele-Shaw
flow is the recent book of Gustafsson and Vasiliev \cite{GV}.
\medskip

\noindent\bf An obstacle problem. \rm
Consider the following function:
$$V(z)=-\weight(z),\qquad z\in\C.$$
For $0<\tau<+\infty$, let ${\mathrm {SH}}_\tau(\C)$ denote the 
collection of (extended real-valued) functions that are superharmonic in 
the whole plane $\C$, and which are of the form
$$-2\tau\,\log|z|+R(z),$$ 
where $R(z)$ is superharmonic in a punctured neighborhood of infinity, and
bounded from below there. If we have two functions 
$f_1,f_2\in{\mathrm {SH}}_\tau(\C)$, then the minimum of the two -- 
$\min\{f_1,f_2\}$ -- is in ${\mathrm {SH}}_\tau(\C)$ as well. In 
addition, the decreasing limit of functions in 
${\mathrm {SH}}_\tau(\C)$ remains in ${\mathrm {SH}}_\tau(\C)$ 
unless it degenerates to $-\infty$ identically. 
{\sl We define $\widehat V_{\tau}$ to be the least majorant 
(or lower envelope) of $V$ in the class ${\mathrm {SH}}_\tau(\C)$.} In other 
words, $V$ is an obstacle, and $\widehat V_{\tau}$ solves the {\sl obstacle 
problem}.
\medskip

\noindent\bf The extremal measure and the obstacle problem. \rm
The next proposition explains the relationship between the obstacle problem
and the solution $\diff\widehat\meas$ to the optimization problem 
(\ref{cap-def}) with $\theta=1/(2\tau)$. For convenience of notation, we shall
write $\widehat\meas_{\weight,\tau}$ in place of 
$\widehat\meas_{\weight,1/(2\tau)}$ when we need to indicate the dependence
on the parameter $\tau$. 

\begin{prop}
The logarithmic potential for the extremal measure $\widehat\meas=
\widehat\meas_{\weight,\tau}$ is
$$L[\widehat\meas](z)=\frac1{2\tau}\,\widehat V_{\tau}(z)+C_\tau,$$
where $C_\tau$ is the constant
$$C_\tau=\frac1{2\tau}\bigg\{\int_\C \weight(z)\,\diff\widehat\meas(z)+
\iint_{\C^2}\log\frac{1}{|z-w|}\,\diff\widehat\meas(z)\diff\widehat\meas(w)
\bigg\}.$$
\label{prop-sigma-env}
\end{prop}

\begin{proof}
This follows from Theorem 4.1 \cite[p. 49]{SaTo}.
\end{proof}

We recall that the function $V=-\Phi$ is assumed to be of class $C^2$ 
throughout $\C$. We shall need the class $C^{1,1}$, which consists of 
continuously differentiable functions whose gradient is locally Lipschitz. 
In other words, all second order partial derivatives of the function are 
locally bounded. 

\begin{prop}
The envelope $\widehat V_\tau$ is of class $C^{1,1}$ throughout $\C$.
\label{prop-C11}
\end{prop}

\begin{proof}
In view of Proposition \ref{prop-sigma-env}, the envelope function 
$\widehat V_\tau$ is harmonic off $\supp\,\widehat\meas$. Since the extremal
measure $\widehat\meas$ has compact support, it follows that we can find an 
open circular disk $\D(0,R)$ of radius $R$, $0<R<+\infty$, around the origin, 
such that $\supp\,\widehat\meas$ is compactly contained in $\D(0,R)$. Along 
the boundary $\T(0,R)$ of the disk, $\widehat V_\tau$ is harmonic, and in 
particular smooth. We consider the following obstacle problem. Suppose $v$
is superharmonic on $\D(0,R)$ and $C^1$-smooth on $\bar\D(0,R)$, and
\begin{equation}
\left\{
\begin{gathered}v(z)\ge V(z),\qquad z\in\D(0,R),\\
v(z)=\widehat V_\tau(z),\qquad z\in\T(0,R).
\end{gathered}
\right.
\label{eq-obstprob}
\end{equation}
The lower envelope of all such $v$ is denoted by $\widetilde V_\tau$. A simple
argument shows that $\widetilde V_\tau$ coincides with $\widehat V_\tau$.
Indeed, we easily see that $\widetilde V_\tau\le\widehat V_\tau$ on $\D(0,R)$,
with equality on $\T(0,R)$; moreover, if we extend $\widetilde V_\tau$ to 
all of $\C$ by declaring that it should equal $\widehat V_\tau$ off 
$\D(0,R)$, we get a superharmonic function, which belongs to 
${\mathrm {SH}}_\tau(\C)$. But -- by definition -- $\widehat V_\tau$ is the
smallest such majorant, and we must have equality: $\widetilde V_\tau=
\widehat V_\tau$. 

In view of the known smoothness properties of solutions of obstacle problems
of the type (\ref{eq-obstprob}), it follows from the $C^2$-smoothness of
$V$ that $\widehat V_\tau$ is of class $C^{1,1}$ on $\D(0,R)$. In the 
rest of the plane, it is harmonic, and therefore automatically of class
$C^{1,1}$.  The basic references on obstacle problems are Chapter 1 in
\cite{Fried}, or the paper \cite{CaKi} by Caffarelli and Kinderlehrer.

The proof is complete.
\end{proof}

We introduce the notation
\begin{equation}
\calD_\Phi(\tau)=\Big\{z\in\C:\,V(z)<\widehat V_\tau(z)\Big\};
\label{eq-HSdom}
\end{equation}
this is an open set which contains a punctured neighborhood of the point
at infinity. 

\begin{prop} 
We have
$$\big\{z\in\C:\,\Delta\weight(z)<0\big\}\subset\calD_\weight(\tau).$$
\end{prop}

\begin{proof}
Let $z_0\in\C$ be a point with $\Delta\weight(z_0)<0$, and suppose that
$z_0\notin\calD_\weight(\tau)$. Then
$$\widehat V(z_0)=V(z_0)=-\weight(z_0),$$
so that if we put $U=\widehat V-V$, we get $U(z_0)=0$, while $U\ge0$
everywhere, and $U$ is strictly superharmonic in a neighborhood of $z_0$.
By the strong maximum principle, then, this is not possible. The result is
immediate.
\end{proof}

We also need the ``harmonicity'' set of $\weight$:
$$\calH_\weight=\big\{z\in\C:\,\Delta\weight(z)=0\big\},$$
which is a closed set in $\C$. 
If $X$ is a Borel measurable subset of the
plane $\C$, let us agree to say that $z_0\in\C$ is a {\sl pseudo-interior
point for} $X$ if there exists a small open disk $D$ around $z_0$ such that
$$|D\cap X|_2=|D|_2$$
holds, where $|\cdot|_2$ is the operation of taking the Lebesgue are measure
of the given set. For instance, all interior points of the set $X$ are 
pseudo-interior points of $X$. 

We are now in a position to reformulate in precise terms our second theorem
from the introduction (Theorem \ref{thm-II}).

\begin{thm}
In the sense of measures, we have
\begin{equation*}
\diff\widehat\meas_{\weight,\tau}(z)=\frac1{4\pi\tau}\,
1_{\C\setminus\calD_\Phi(\tau)}(z)\,
\Delta\weight(z)\,\diff\lambda_2(z),\qquad z\in\C.
\end{equation*}
In particular, the support of the measure $\widehat\meas_{\weight,\tau}$ 
consists of all points of $\C$ which are not pseudo-interior for 
the set $\calD_\weight(\tau)\cup\calH_\weight$.
\label{prop-basic}
\end{thm}

\begin{proof}
We get from Proposition \ref{prop-sigma-env} that
\begin{equation}
-2\pi\,\diff\widehat\meas(z)=\Delta L[\widehat\meas](z)\,\diff\lambda_2(z)
=\frac1{2\tau}\,\Delta\widehat V_{\tau}(z)\,\diff\lambda_2(z).
\label{eq-aa1}
\end{equation}
According to \cite[p. 53]{KiS}, we have 
\begin{equation}
\Delta\widehat V_\tau(z)=\Delta V(z)=-\Delta\weight(z),\qquad 
z\in\C\setminus\calD_\weight(\tau),
\label{eq-aa2}
\end{equation}
in the almost-everywhere sense. On the non-coincidence set, however, 
$\widehat V_\tau$ is harmonic:
\begin{equation}
\Delta\widehat V_\tau(z)=0,\qquad z\in\calD_\weight(\tau).
\label{eq-aa3}
\end{equation}
This follows by a standard Perron process argument (see, for instance, 
Proposition 2.2 in \cite{HS}). By Proposition \ref{prop-C11}, the function 
$\Delta\widehat V_\tau$ is in $L^\infty_{\mathrm loc}(\C)$, so that we see
from (\ref{eq-aa1}), (\ref{eq-aa2}), and (\ref{eq-aa3}) that  
\begin{equation*}
\diff\widehat\meas_{\weight,\tau}(z)=\frac1{4\pi\tau}\,
1_{\C\setminus\calD_\Phi(\tau)}(z)\,
\Delta\weight(z),\qquad z\in\C,
\end{equation*}
as asserted. The statement regarding the support of 
$\widehat\meas_{\weight,\tau}$ is an easy consequence of this identity.
The proof is complete.
\end{proof}

The obstacle problem considered here share many features in common with
the Hele-Shaw flow domains considered by Hedenmalm and Shimorin in \cite{HS}
However, the time parameter $\tau$ flows in the opposite direction as compared
with the time parameter $t$ used in \cite{HS}. This means that any cusps
of the compact set $\C\setminus\calD_\weight(\tau)$ should point outward. 
Also, if $\weight$ is real-analytic, the smoothness analysis of based on
Sakai's work \cite{Sakai} done in \cite{HS} carries over to the situation 
treated here. 

The domains $\calD(\tau)$ grow as $\tau$ {\sl decreases}, and the 
``harmonic moments'' are preserved, with the exception of the first;
this is the content of the following proposition. 

\begin{prop}
Let $0<\tau<\tau'<+\infty$. Suppose $h$ is harmonic and bounded in 
$\calD_\weight(\tau)$, with an extension to $\C$ that is locally of Sobolev
class $W^{2,2}$. We then have the equality
$$4\pi\,(\tau'-\tau)\,h(\infty)=\int_{\calD_\weight(\tau)\setminus
\calD_\weight(\tau')}h(z)\,\Delta\weight(z)\,\diff\lambda_2(z).$$
\label{prop-III}
\end{prop}

\begin{proof}
This follows from a suitable application of Green's formula, analogous
to what is done in \cite{HS}.
\end{proof}

If the domains $\calD_\weight(\tau)$ are smooth and vary smoothly with 
$\tau$, then, by arguing as in \cite{HS}, pp. 188-189, we see that the 
assertion of Proposition \ref{prop-III} has the interpretation that the 
boundary $\partial\calD_\weight(\tau)$ propagates with velocity proportional 
to a weight times the gradient of the Green function with singularity at 
infinity. 

The growth of the complementary sets $\C\setminus\calD_\weight(\tau)$
as $\tau$ increases is quite an interesting process. For $\tau$ close
to $0$, we should expect each such complementary set to be localized
close to the points where the minimum of the potential function $\weight$
is attained. If there is only one minimum point $z_0$, and the function 
$\weight$ is convex in a neighborhood of $z_0$, then we can show that
a small neighborhood of $z_0$ is contained in 
$\C\setminus\calD_\weight(\tau)$, for each fixed $\tau$, by comparing with
a concave majorant. An interesting question seems to be whether
$$\bigcap_{0<\tau<+\infty}\calD_\weight(\tau)=\emptyset.$$
It is easy to see that a necessary condition for this to happen is that
$\weight$ be subharmonic everywhere. We do not know to what extent this 
condition is sufficient. 

\section{Condensation of quantum Hele-Shaw flow}

\noindent \bf The continuous limit of eigenvalues of normal matrices. \rm
For positive real $\beta$, let $\mu_{\beta,\weight}$ denote the measure
\begin{equation}
\diff\mu_{\beta,\weight}(z)=
\exp\big\{-\beta\,\weight(z)\big\}\,\diff\lambda_2(z),
\qquad z\in\C,
\label{eq:qs-1}
\end{equation}
where we recall that $\diff\lambda_2$ is area measure in the plane. Throughout
this section, we shall use -- as before -- the following scaling choice of 
$\beta$ as we increase $N$:
\begin{equation}
\beta=2\theta\,(N-1)=\frac{N-1}\tau,
\label{eq-beta}
\end{equation}
where $\tau$ is a fixed positive real parameter.   
Also, we recall the previously used notation 
$$\zed_{[1,n]}=(z_1,\ldots,z_n),\qquad \zed_{[1,N]}=(z_1,\ldots,z_N).$$
Given this setup, we consider the $n$-point correlation measure 
$\corr^{(n)}_{\mu,N}$ with $\mu=\mu_{\beta,\weight}$ and $\beta$ given by 
(\ref{eq-beta}). We want to understand the asymptotic behavior of this 
measure as $N\to+\infty$. This means that we should understand the impact 
of working with the $L^2$ norm in place of the $L^\infty$ norm in the 
setting of (\ref{distr-pi-3}).
\medskip

We recall that $\weight$ is a $C^2$-smooth real-valued function, which tends
to infinity at infinity at a pace prescribed by (\ref{est-below}). 
It follows from this that there exists a compact subset $K$ of $\C$ such that
\begin{multline}
\kappa_{\weight}(\theta)+1<
\log\frac1{E_{\weight}(z,w;\theta)}\\
=\theta\,\big[\weight(z)+\weight(w)\big]
+\log\frac1{|z-w|},
\qquad (z,w)\in\C^2\setminus(K\times K),
\label{def-K}
\end{multline}
where we recall the definition (\ref{cap-def}) of $\kappa_{\weight}(\theta)$.
In analogy with the energy function (\ref{energy}), we introduce the function
$\calE^\sharp_\weight\big[\zed_{[1,N]}\big](\theta)$,
\begin{multline}
\calE^\sharp_\weight[\zed_{[1,N]}](\theta)=
\frac2{N(N-1)}\,\log\prod_{j,k:\,j<k}\frac1{E_\weight(z_j,z_k;\theta)}\\
=\frac1{N(N-1)}\sum_{j,k:\,j\neq k}\log\frac1{E_\weight(z_j,z_k;\theta)},
\label{energy-sharp}
\end{multline} 
where $E_\weight(z,w;\theta)$ is given by (\ref{E-fn}), and 
$\theta=1/(2\tau)$. We recall that 
\begin{equation}
\log\frac1{M_{\weight,N}(\theta)}\le
\calE^\sharp_\weight[\zed_{[1,N]}](\theta),\qquad \zed_{[1,N]}\in\C^N,
\label{eq-estenergy}
\end{equation}
with equality only at the weighted Fekete points. 
For positive real $\epsilon$, let
$$\Aset_{\weight,N}(\epsilon,\theta)=\bigg\{\zed_{[1,N]}\in\C^N:\,\,
\calE^\sharp_\weight[\zed_{[1,N]}](\theta)\le
\log\frac1{M_{\weight,N}(\theta)}+\epsilon \bigg\}.$$
We need to know that the proportion of points in 
$\Aset_{\weight,N}(\epsilon,\theta)$ which stays in the compact set 
$K$ converges to $1$ as $N\to+\infty$ and $\epsilon\to0$. 

\begin{prop}
Suppose $\zed_{[1,N]}\in\Aset_N(\epsilon,\theta)$, and that 
$0<\epsilon<\frac12$. Let $N_K$ denote the number of indices $j$ for which 
$z_j\in K$. Then, for sufficiently large $N$, we have
$$\frac{N_K}{N}\ge1-2\epsilon.$$
\label{prop-eps}
\end{prop}

\begin{proof}
Let $X\subset\{1,\ldots,N\}$ be the subset of all indices $j$ for which
$z_j\in K$, and let $Y$ be the complement in $\{1,\ldots,N\}$. We split
the sum defining $\calE^\sharp_\weight\big[\zed_{[1,N]}\big](\theta)$
accordingly:
\begin{multline*}
\calE^\sharp_\weight\big[\zed_{[1,N]}\big](\theta)
=\frac1{N(N-1)}
\sum_{j,k\in X:\,j\neq k}\log\frac1{E_\weight(z_j,z_k;\theta)}\\
+\frac2{N(N-1)}\sum_{j\in X,\,k\in Y:\,j\neq k}
\log\frac1{E_\weight(z_j,z_k;\theta)}
\\
+\frac1{N(N-1)}
\sum_{j,k\in Y:\,j\neq k}\log\frac1{E_\weight(z_j,z_k;\theta)}.
\end{multline*}
The sum of the last two terms is estimated from below by
$$\bigg(1-\frac{N_K(N_K-1)}{N(N-1)}\bigg)\big(\kappa_\weight(\theta)+1\big),$$
while the first term of the right hand side is estimated in the following
manner:
\begin{multline*}
\log\frac1{M_{\weight,N_K}(\theta)}\le
\calE^\sharp_\weight[\zed_{[1,N]}](\theta)\\
=\frac{1}{N_K(N_K-1)}
\sum_{j,k\in X:\,j\neq k}\log\frac1{E_\weight(z_j,z_k;\theta)},
\qquad \zed_{[1,N]}\in\C^N.
\end{multline*}
We now see that for $\zed_{[1,N]}\in\Aset_N(\epsilon,\theta)$,
\begin{multline*}
\frac{N_K(N_K-1)}{N(N-1)}
\log\frac1{M_{\weight,N_K}(\theta)}+
\bigg(1-\frac{N_K(N_K-1)}{N(N-1)}\bigg)\big(\kappa_\weight(\theta)+1\big)\\
\le
\calE^\sharp_\weight[\zed_{[1,N]}](\theta)\le 
\log\frac1{M_{\weight,N}(\theta)}+\epsilon.
\end{multline*}
It follows that
\begin{equation}
1+\kappa_\weight(\theta)-\log\frac1{M_{\weight,N}(\theta)}-\epsilon
\le
\frac{N_K(N_K-1)}{N(N-1)}\left(1+\kappa_\weight(\theta)-
\log\frac1{M_{\weight,N_K}(\theta)}\right).
\label{eq-est}
\end{equation}
For large values of $N$, the left hand side is positive, and hence the right
hand side is, too. We quickly rule out the possibility that the first on the
right hand side is negative (which corresponds to $N_K=0$), 
in which case we are left with both factors on the right hand side being 
positive. We get from (\ref{eq-est}) that
\begin{equation}
\frac{N_K(N_K-1)}{N(N-1)}\ge
1-\frac{\epsilon}{1+\kappa_\weight(\theta)-
\log\frac1{M_{\weight,N_K}(\theta)}},
\label{eq-est2}
\end{equation}
if we recall that $M_{\weight,N}(\theta)\le M_{\weight,N_K}(\theta)$.
Next, we note that $N_K\to+\infty$ as $N$ tends to infinity; therefore, for
sufficiently large $N$, 
$$0\le\log\frac1{M_{\weight,N_K}(\theta)}-\kappa_\weight(\theta)\le\frac12.$$ 
The assertion of the proposition now follows from (\ref{eq-est2}).
\end{proof}

\begin{lem}
Let $\meas\in\calP_c(\C)$ be absolutely continuous
with respect to two-dimensional Lebesgue measure, 
$$\diff\meas(z)=S(z)\,\diff\lambda_2(z),$$
where $S\ge0$ is area-summable and $S\log^+ S$ is area-summable as well.
The normalization constant $Z_N$ in the definition of 
$\diff\Pi_{N,\beta,\weight}(\zed_{[1,N]})$ then has the following bound:
\begin{equation*}
Z_N\ge\exp \bigg\{-N(N-1)\,\calE_\weight[\meas](\theta)
-N\int_\C \log S(z)\,\diff\meas(z)\bigg\},
\end{equation*}
with the understanding that $\log S\, \diff\meas$ vanishes off the
support of $\meas$. 
\label{lem-1.1}
\end{lem}

\begin{proof}
We recall first the definition of $Z_N$:
\begin{equation}
Z_N=\int_{\C^N}\exp\bigg\{-\beta\sum_{j=1}^N \weight(z_j)\bigg\}\,
\big|\triangle(\zed_{[1,N]})\big|^2\,\diff \lambda_{2N}(\zed_{[1,N]}).
\label{eq-Z}
\end{equation}
Let $\Meas$ denote the set where $S(z)>0$; as $\meas$ has
compact support, we may assume that the set $\Meas$ is bounded. 
We then have
\begin{multline*}
Z_N\ge\int_{\Meas^N}
\exp\bigg\{-\beta\sum_{j=1}^N \weight(z_j)\bigg\}\,
\big|\triangle(\zed_{[1,N]})\big|^2\,\frac1{S(z_1)\cdots S(z_N)}\,
\diff\meas(z_1)\cdots\diff\meas(z_N)\\
\ge\exp \Bigg\{\int_{\Meas^N}
\bigg(-\beta\sum_{j=1}^N \weight(z_j)+
2\log\big|\triangle(\zed_{[1,N]})\big|\,-\sum_{j=1}^N
\log S(z_j)\bigg)\\
\times\diff\meas(z_1)\cdots\diff\meas(z_N)
\Bigg\},
\end{multline*}
where the second estimate is due to Jensen's inequality. In view of
(\ref{eq-beta}) and the definition of the energy functional, we have
\begin{equation*}
Z_N\ge\exp \bigg\{-N(N-1)\,\calE_\weight[\meas](\theta)
-N\int_\Meas \log S(z)\,\diff\meas(z)\bigg\},
\end{equation*}
as claimed.
\end{proof}

\begin{rem}
We note that in view of Theorem \ref{prop-basic}, the extremal measure
$\meas=\widehat\meas$ enjoys the regularity property of Lemma 
\ref{lem-1.1}. 
\label{rem-1.0}
\end{rem}

The next proposition is central to the argument; its proofs mimics
Johansson's method for eigenvalues of self-adjoint matrices \cite{Kurt}. 

\begin{prop} For sufficiently large $N$, we have 
$$\Pi_{N,\beta,\weight}\big(\C^N\setminus\Aset_N(\epsilon,\theta)\big)\le
e^{-\epsilon\,N(N-1)/2}.$$
\label{prop-picompl}
\end{prop}

\begin{proof}
For 
$\zed_{[1,N]}\in\C^N\setminus\Aset_N(\epsilon,\theta)$,
we have 
\begin{multline*}
\calE^\sharp_\weight\big[\zed_{[1,N]}\big](\theta)=
\frac1{N(N-1)}\sum_{j,k:\,j\neq k}
\log\frac1{E_\weight(z_j,z_k;\theta)}
>\log\frac1{M_{\weight,N}(\theta)}+\epsilon\\
\ge \kappa_\weight(\theta)+
\epsilon.
\end{multline*}
We recall the estimate (\ref{est-below}) from below of the weight $\weight$,
which for an appropriate choice of the parameter $A$ gives
\begin{equation*}
\calE^\sharp_\weight\big[\zed_{[1,N]}\big](\theta)\ge-C+\frac1{N}
\sum_{j=1}^N\log\big(1+|z_j|^2\big), 
\qquad \zed_{[1,N]}\in\C^N,
\end{equation*}
for a positive constant $C$ that only depends on $\theta$ and $\weight$. 
We form a convex combination of these two estimates ($0\le\gamma\le1$):
\begin{multline}
\calE^\sharp_\weight\big[\zed_{[1,N]}\big](\theta)\ge
(1-\gamma)\big[\kappa_{\weight}(\theta)+\epsilon\big]-C\gamma\\
+\frac{\gamma}{N}\sum_{j=1}^N\log\big(1+|z_j|^2\big),\qquad
\zed_{[1,N]}\in\C^N\setminus\Aset_N(\epsilon,\theta).
\label{eq-a2}
\end{multline}
We recall the formula defining $\Pi_{N,\beta,\weight}$, while keeping 
in mind the identity (\ref{eq-dens}),
\begin{multline}
\diff\Pi_{N,\beta,\weight}\big(\zed_{[1,N]}\big)=\frac1{Z_N}\bigg[
\prod_{j,k:\,j\neq k} E_\weight(z_j,z_k;\theta)\bigg]\,\diff\lambda_{2N}
\big(\zed_{[1,N]}\big)\\
=\frac1{Z_N}\,\exp\Big\{-N(N-1)\,\calE^\sharp_\weight
\big[\zed_{[1,N]}\big](\theta)\Big\}\,
\diff\lambda_{2N}\big(\zed_{[1,N]}\big).
\label{eq-pi3}
\end{multline}
So, on the set $\C^N\setminus\Aset_N(\epsilon,\theta)$,
we get, in view of (\ref{eq-a2}),
\begin{multline*}
\diff\Pi_{N,\beta,\weight}\big(\zed_{[1,N]}\big)\le\frac1{Z_N}\,
\exp\bigg\{-N(N-1)(1-\gamma)\big[\kappa_{\weight}(\theta)+\epsilon\big]\\
+CN(N-1)\gamma-\gamma(N-1)\sum_{j=1}^N\log\big(1+|z_j|^2\big)
\bigg\}\,
\diff\lambda_{2N}\big(\zed_{[1,N]}\big). 
\end{multline*}
By Lemma \ref{lem-1.1}, we have
\begin{equation*}
Z_N\ge\exp \bigg\{-N(N-1)\,\kappa_\weight(\theta)
-B\,N\bigg\},
\end{equation*}
where $B$ is the real number
$$B=\int_\C \log\Delta\weight(z)\,\diff\widehat\meas(z)
+\log\frac{\theta}{\pi}.$$
As we combine this with (\ref{eq-a2}), we arrive at
\begin{multline*}
\diff\Pi_{N,\beta,\weight}\big(\zed_{[1,N]}\big)\le
\exp\bigg\{-N(N-1)\big[\gamma\,\kappa_{\weight}(\theta)+
(1-\gamma)\epsilon\big]+B\,N\\
+C\gamma N(N-1)-\gamma(N-1)\sum_{j=1}^N\log\big(1+|z_j|^2\big)
\bigg\}\,
\diff\lambda_{2N}\big(\zed_{[1,N]}\big)\\
=\exp\bigg\{-N(N-1)\big[\gamma\,\kappa_{\weight}(\theta)+
(1-\gamma)\epsilon\big]+B\,N+C\gamma N(N-1)\bigg\}\\
\times\prod_{j=1}^N\big(1+|z_j|^2\big)^{-\gamma(N-1)}
\diff\lambda_{2N}\big(\zed_{[1,N]}\big)
\end{multline*}
on the set $\C^N\setminus\Aset_N(\epsilon,\theta)$. Considering
that
$$\int_\C\big(1+|z|^2\big)^{-\gamma(N-1)}\diff\lambda_2(z)=
\frac{\Pi}{\gamma(N-1)-1},$$
provided that $1/(N-1)<\gamma\le1$, we find that
\begin{multline*}
\Pi_{N,\beta,\weight}\big(\C^N\setminus\Aset_N(\epsilon,\theta)\big)\le
\exp\bigg\{-N(N-1)\big[\gamma\,\kappa_{\weight}(\theta)+
(1-\gamma)\epsilon-C\gamma\big]+B\,N\bigg\}
\\
\times\bigg(\frac{\Pi}{\gamma(N-1)-1}\bigg)^N.
\end{multline*}
The dominant contribution in the expression which is exponentiated is
$$-N(N-1)\big[\gamma\,\kappa_{\weight}(\theta)+
(1-\gamma)\epsilon-C\gamma\big].$$
We would like to pick $\gamma$, $1/(N-1)<\gamma<1$, such that
$$0<\gamma\,\kappa_{\weight}(\theta)+
(1-\gamma)\epsilon-C\gamma.$$
This is possible; indeed, without loss of generality, we may assume that
$C$ is greater than $\kappa_\weight(\theta)$, in which case one such choice 
is 
$$\gamma=\frac{\epsilon}{2(C-\kappa_\weight(\theta)+\epsilon)}.$$
This value of $\gamma$ yields the assertion of the proposition.
The proof is complete.
\end{proof}

For a point $\zed_{[1,N]}\in\C^N$, we define the associated weighted
sum of point masses $\meas[\zed_{[1,N]}]\in\calP_c(\C)$ by the 
formula
\begin{equation}
\diff\meas[\zed_{[1,N]}](z)=\frac1{N}\sum_{j=1}^N\diff\delta_{z_j}(z),
\qquad z\in\C,
\label{eq-sigmaN}
\end{equation}
where $\delta_w$ means the Dirac point mass at $w\in\C$. Also, let
$C_b(\C)=C(\C)\cap L^\infty(\C)$ denote the space bounded complex-valued
continuous functions on $\C$. 

\begin{prop}
Suppose $\meas_N=\meas\big[\zed_{[1,N]}\big]$ is as above, with 
$$\zed_{[1,N]}=\big(z_1,\ldots,z_N\big)\in\C^N.$$
Suppose, moreover, that 
$$\calE^\sharp_\weight\big[\zed_{[1,N]}\big](\theta)
\to\kappa_\weight(\theta)$$
as $N\to+\infty$. Then $\meas_N$ converges to $\widehat\meas$, the 
extremal measure, in the weak-star topology, as $N\to+\infty$. 
In other words, for each $f\in C_b(\C)$, we have
$$\int_\C f(z)\,\diff\meas_N(z)\to \int_\C f(z)\,\diff\widehat\meas(z)\quad
\text{as}\quad N\to+\infty.$$
\label{prop-limitmeas}
\end{prop}

\begin{proof}
First, fix a small but positive $\epsilon$. By assumption, for sufficiently 
big $N$,
$$\zed_{[1,N]}\in\Aset_N(\epsilon,\theta),$$
so that by Proposition \ref{prop-eps}, $\meas_N(K)\ge1-2\epsilon$. As
$\epsilon$ can be made as small as we like, and each $\meas_N$ is a 
probability measure, we see that
$$\meas_N(K)\to1\quad\text{as}\quad N\to+\infty.$$

The space of all finite complex Borel measures on $K$ is a Banach space, with
weak-star compact unit ball. This means that
each subsequence of the sequence $\meas_1|_K,\meas_2|_K,{\meas_3}|_K,
\ldots$ has a weak-star convergent subsequence. We shall show that any such 
limit coincides with $\widehat\meas$, from which the assertion follows.

Without loss of generality, then, by passing to a subsequence, we may assume 
that the sequence $\meas_N|_K$, $N=1,2,3,\ldots$, converges weak-star
itself to a limit, which we call $\widetilde\meas$. The weak-star limit of 
$\meas_N$ is then also $\widetilde\meas$, given that 
$\meas_N(\C\setminus K)\to0$ as $N\to+\infty$.
By testing with the function $f=1$, we see that $\widetilde\meas$ is
a probability measure supported on $K$, so that
$\widetilde\meas\in\calP_c(\C)$. {\sl We claim that} 
$\calE_\weight\big[\widetilde\meas\big](\theta)\le\kappa_\weight(\theta)$.
Once this is established, the equality $\widetilde\meas=\widehat\meas$
follows from the uniqueness of the extremal measure. 
We use a cut-off argument: we note that
\begin{equation*}
\calE^\sharp_\weight\big[\zed_{[1,N]}\big](\theta)=\frac{2\theta}{N}
\sum_{j=1}^N\weight(z_j)+\frac1{N(N-1)}\sum_{j,k:\,j\neq k}
\log\frac{1}{|z_j-z_k|},
\end{equation*}
so that if $L$ is a real parameter,
\begin{multline*}
\calE^\sharp_\weight\big[\zed_{[1,N]}\big](\theta)\ge
\frac{2\theta}{N}
\sum_{j=1}^N\weight(z_j)
+\frac1{N(N-1)}\sum_{j,k=1}^N
\min\bigg\{\log\frac{1}{|z_j-z_k|},L\bigg\}\\
-\frac{L}{N-1}.
\end{multline*} 
Now, let $N\to+\infty$, so that we get, in view of our assumptions, 
\begin{equation*}
\kappa_\weight(\theta)\ge2\theta\int_\C\weight(z)\,\diff\widetilde\meas(z)
+\int_{\C^2}\min\bigg\{\log\frac{1}{|z-w|},L\bigg\}\,\diff\widetilde\meas(z)
\,\diff\widetilde\meas(w).
\end{equation*}
As we let $L$ tend to $+\infty$, the right hand side approaches 
$\calE_\weight\big[\widetilde\meas\big](\theta)$, whence the claim is
immediate.
\end{proof}

As in the introduction, we let $\Pi^{(n)}_{N,\beta,\weight}$ be the
$n$-point marginal distribution measure and $\corr^{(n)}_{N,\beta,\weight}$ 
the $n$-point correlation measure for the probability distribution measure
$\Pi_{N,\beta,\weight}$. We keep $\beta$ connected with $\tau$, 
$N$, and $\theta$ via (\ref{eq-beta}), and fix $\theta$ (or, if you like, 
$\tau$) while $N$ grows.
\medskip

Let $C_b(\C^n)=C(\C^n)\cap L^\infty(\C^n)$ denote the Banach space of bounded
complex-valued continuous functions on $\C^n$. We arrive at a precise 
reformulation of our first main result (Theorem \ref{thm-I}).

\begin{thm} We have
$$\Pi^{(n)}_{N,\beta,\weight}(\C^n)=1,$$
while, as $N\to+\infty$,
$$\Pi^{(n)}_{N,\beta,\weight}(K^n)\to 1.$$ 
Moreover, for each $f\in C_b(\C^n)$, we have, as $N\to+\infty$,
\begin{equation}
\int_{\C^n}f\big(\zed_{[1,n]}\big)\,
\diff\Pi^{(n)}_{N,\beta,\weight}\big(\zed_{[1,n]}\big)
\to\int_{\C^n}f\big(\zed_{[1,n]}\big)\,
\diff\widehat\meas(z_1)\cdots\diff\widehat\meas(z_n).
\label{eq-limitthm}
\end{equation}
In other words, in the weak-star topology of measures, we have, as 
$N\to+\infty$,
$$\diff\Pi^{(n)}_{N,\beta,\weight}\big(\zed_{[1,n]}\big)\to
\diff\widehat\meas(z_1)\cdots\diff\widehat\meas(z_n).$$ 
Here, $\widehat\meas=\widehat\meas_{\weight,\theta}$ is the extremal 
measure.
\label{thm-main}
\end{thm}

\begin{proof}
We note that the total mass of the measure $\Pi^{(n)}_{N,\beta,\weight}$
is $1$. 

As a second step, we prove (\ref{eq-limitthm}) under the slightly more 
restrictive assumption $f\in C_c(\C)$, which means that the test function 
$f$ has compact support. 

Let $\varsigma$ be a permutation of $\{1,2,3,\ldots,N\}$. Then, due to the
symmetry properties of $\Pi_{N,\beta,\weight}$, 
\begin{multline*}
\int_{\C^n}f\big(\zed_{[1,n]}\big)\,
\diff\Pi^{(n)}_{N,\beta,\weight}\big(\zed_{[1,n]}\big)=
\int_{\C^N}f(z_1,\ldots,z_n)\,
\diff\Pi_{N,\beta,\weight}\big(\zed_{[1,N]}\big)\\
=\int_{\C^N}f(z_{\varsigma(1)},\ldots,z_{\varsigma(n)})\,
\diff\Pi_{N,\beta,\weight}\big(\zed_{[1,N]}\big),
\end{multline*}
from which we quickly deduce that
\begin{equation*}
\int_{\C^n}f\big(\zed_{[1,n]}\big)\,
\diff\Pi^{(n)}_{N,\beta,\weight}\big(\zed_{[1,n]}\big)
=\frac1{N!}
\sum_\varsigma \int_{\C^N}
f\big(z_{\varsigma(1)},\ldots,z_{\varsigma(n)}\big)\,
\diff\Pi_{N,\beta,\weight}\big(\zed_{[1,N]}\big),
\end{equation*}
where the sum runs over all permutations of $\{1,2,3,\ldots,N\}$. We
fix a small positive $\epsilon$, and split the integral:
\begin{multline}
\int_{\C^n}f\big(\zed_{[1,n]}\big)\,
\diff\Pi^{(n)}_{N,\beta,\weight}\big(\zed_{[1,n]}\big)\\
=\frac{1}{N!}\int_{\Aset_N(\epsilon,\theta)}\sum_\varsigma 
f\big(z_{\varsigma(1)},\ldots,z_{\varsigma(n)}\big)\,
\diff\Pi_{N,\beta,\weight}\big(\zed_{[1,N]}\big)\\
+\frac{1}{N!}\sum_\varsigma
\int_{{\C^N\setminus\mathfrak A}_N(\epsilon,\theta)}
f\big(z_{\varsigma(1)},\ldots,z_{\varsigma(n)}\big)\,
\diff\Pi_{N,\beta,\weight}\big(\zed_{[1,N]}\big).
\label{eq-splitint}
\end{multline}
By Proposition \ref{prop-picompl}, the last term is bounded in modulus by
\begin{equation}
e^{-\epsilon\,N(N-1)/2}\,\|f\|_{L^\infty(\C)},
\label{eq-est3}
\end{equation}
for large $N$. To deal with the first term on the right hand side, 
we should understand the behavior of
\begin{equation}
\frac{1}{N!}
\sum_\varsigma f\big(z_{\varsigma(1)},\ldots,z_{\varsigma(n)}\big),
\qquad \zed_{[1,N]}\in\Aset_N(\epsilon,\theta).
\label{eq-sum}
\end{equation}
Let us consider the simplest case $n=1$ first. Then (\ref{eq-sum}) amounts
to 
\begin{equation}
\frac1{N}\sum_{j=1}^N f\big(z_{j}\big),
\qquad \zed_{[1,N]}\in\Aset_N(\epsilon,\theta).
\label{eq-sum2}
\end{equation}
By letting $\epsilon$ approach $0$ slowly as $N\to+\infty$, we may
ensure that (\ref{eq-est3}) tends to $0$ as $N\to+\infty$, while 
(\ref{eq-sum2}) approaches
$$\int_\C f(z)\,\diff\widehat\meas(z).$$
It is easy to check that the latter statement entails
$$\frac{1}{N!}\int_{\Aset_N(\epsilon,\theta)}\sum_{j=1}^N 
f\big(z_{j}\big)\,\diff\Pi_{N,\beta,\weight}\big(\zed_{[1,N]}\big)\to
\int_\C f(z)\,\diff\widehat\meas(z)$$
as $N\to+\infty$, if we allow $\epsilon$ to approach $0$ slowly.

The remaining case $n>1$ is handled in an analogous manner: Proposition
\ref{prop-limitmeas} should be replaced by a multidimensional analogue,
which we get by iterated integration. It is useful to keep in mind that 
for large $N$ and fixed $n$, the following collections of sequences are 
asymptotically (as $N\to+\infty$) the same: 
$$\big(\varsigma(1),\ldots,\varsigma(n)\big)\quad\text{and}\quad
\big(j(1),\ldots,j(n)\big),$$
where $\varsigma$ runs over all permutations of $\{1,\ldots,N\}$, and
$j$ runs over all functions $j:\{1,\ldots,n\}\to\{1,\ldots,N\}$. For instance,
$$\frac{N^n(N-n)!}{N!}\to1\quad\text{as}\quad N\to+\infty.$$

There are a couple of assertion that remain to check. If we notice that
the support of $\widehat\meas$ is actually contained in the interior of the
compact set $K$, then, by choosing a smooth cut-off function 
$f=\chi\in C_c(\C^n)$ with $0\le \chi(\zed_{[1,n]})\le1$ everywhere, 
$\chi(\zed_{[1,n]})=1$ on $(\supp\widehat\meas)^n$,
and $\chi(z)=0$ off $K^n$, we find from (\ref{eq-limitthm}) that
$$\liminf_{N\to+\infty}\,\,\Pi^{(n)}_{N,\beta,\weight}(K^n)\ge1.$$ 
Since we are dealing with probability measures, it follows that 
$$\lim_{N\to+\infty}\Pi^{(n)}_{N,\beta,\weight}(K^n)=1,$$  
as claimed. 

It remains to verify that (\ref{eq-limitthm}) holds for all $f\in C_b
(\C^n)$. To this end, we use a smooth cut-off function $\chi$ similar to 
what we defined above, and write $f=\chi\, f+(1-\chi)\,f$. Then 
$\chi\,f\in C_c(\C)$, while the integral of the remainder $(1-\chi)\,f$
tends to $0$ as $N\to+\infty$, by a simple estimate. 
 
The proof is complete.
\end{proof}

\begin{rem}
Let us think of the $N$-tuple $\zed_{[1,N]}=(z_1,\ldots,z_N)$ as a random 
variable taking values in $\C^N$, with distribution measure 
$\Pi_{N,\beta,\weight}$. The marginal distribution of the first $n$ 
($1\le n\le N$) coordinates $\zed_{[1,n]}=(z_1,\ldots,z_n)$ is then given by 
the measure $\Pi_{N,\beta,\weight}^{(n)}$; to stress the dependence on $N$, 
let us write 
$$\zed_{[1,n]}|_N=\big(z_1|_N,\ldots,z_n|_N\big)$$
for the marginal random variable. Now, if $\beta$ grows with $N$ according to
(\ref{eq-beta}), for some fixed positive $\theta$, Theorem 
\ref{thm-main} says that as $N\to+\infty$, $z_1|_N,\ldots,z_n|_N$ become 
identically distributed independent random variables, each having 
$\widehat\sigma_{\weight,\theta}$ as distribution measure.  In particular, 
the covariance type matrix (with $\mu=\mu_{\beta,\weight}$)
$$\,{\mathbf K}^n_{\mu,N}(\zed_{[1,n]},\zed_{[1,n]})=   
\Big(\,K_{\mu,N}(z_j,z_k)\Big)_{j,k=1}^n$$
becomes asymptotically diagonal as $N\to+\infty$.  
\end{rem}

\begin{ex} We need to give an example to illustrate the result. We consider 
the rather trivial case of $\weight(z)=|z|^2$, and put $n=1$. Then 
$\Delta\weight(z)=4$ is constant. A computation
shows that with $\mu=\mu_{\beta,\Phi}$ given by (\ref{eq-mu}), the 
reproducing kernel function for the polynomial subspace is 
$$K_{\mu,N}(z,w)=
\frac{\beta}{\pi}\sum_{j=0}^{N-1}\frac{(\beta z\bar w)^j}{j!}.$$
We then have
$$\diff\Pi^{(n)}_{N,\beta,\weight}(z)=\frac1{N}\,K_\mu(z,z)\,e^{-\beta|z|^2}
\,\diff\lambda_2(z)=\frac{\beta}{N\pi}
\,e^{-\beta|z|^2}\sum_{j=0}^{N-1}\frac{\big(\beta |z|^2\big)^j}{j!}
\,\diff\lambda_2(z),$$
while it is easy to check that 
$$\calD_\weight(\tau)=\big\{z\in\C:\,|z|>\sqrt{\tau}\big\}.$$
In view of this, the content of Theorem \ref{thm-main} in this simple case is 
that, for fixed positive $\tau$, the function
$$e^{-(N-1)|z|^2/\tau}\sum_{j=0}^{N-1}\frac1{j!}\,
\left(\frac{(N-1)|z|^2}{\tau}\right)^j,$$
which is real-valued, with values between $0$ and $1$, tends to $1$ as
$N\to+\infty$ for $|z|<\sqrt{\tau}$, and to $0$ for $|z|>\sqrt{\tau}$. 
This fact is of course well-known. In \cite{ElFe}, the domain 
$\calD_\weight(\tau)$ is computed explicitly and shown to be an ellipse
in the more general case (still with constant $\Delta\Phi(z)=4$)
$$\Phi(z)=|z|^2+a\re(z^2),\qquad -1<a<1.$$

\end{ex}

\bigskip

\noindent \bf Acknowledgements. \rm The authors wish to thank Bj\"orn 
Gustafsson, Kurt Johansson, Sergei Naboko, and Henrik Shahgholian for helpful 
conversations while this work progressed.

\bigskip


\bigskip

\begin{thebibliography}{KKK}

\bibitem{ABWZ} O. Agam, E. Bettelheim, P. Wiegmann, A. Zabrodin, 
{\it Viscous fingering and the shape of an electronic droplet 
in the Quantum Hall regime}, Phys. Rev. Lett. {\bf88} (2002), 236801(1--4).

\bibitem{ARS} A. Aleman, S. Richter, C. Sundberg, \it The majorization 
function and the index of invariant subspaces in the Bergman spaces.  \rm
J. Anal. Math. {\bf 86} (2002), 139--182. 

\bibitem{CaKi} L. A. Caffarelli, D. Kinderlehrer, \it Potential methods in
variational inequalities, \rm J. Anal. Math. {\bf37} (1980), 285--295.

\bibitem{ChZa} L.-L. Chau, O. Zaboronsky, \it On the structure of
correlation functions in the normal matrix model, 
\rm Comm. Math. Phys. {\bf196} (1998), no. 1, 203--247. 

\bibitem{Deift} P. Deift, \it Orthogonal polynomials and random matrices: 
a Riemann-Hilbert approach. \rm Courant Lecture Notes in Mathematics, {\bf3}. 
New York University, Courant Institute of Mathematical Sciences, New York; 
American Mathematical Society, Providence, RI, 1999.

\bibitem{ElFe} P. Elbau, G. Felder, \it Density of eigenvalues of random
normal matrices. \rm Preprint 2004.

\bibitem{Ful} W. Fulton, \it Eigenvalues, invariant factors, highest weights, 
and Schubert calculus. \rm  Bull. Amer. Math. Soc. (N.S.) {\bf37} (2000), 
no. 3, 209--249.

\bibitem{Fried} A. Friedman, \it Variational principles and free-boundary
problems, \rm 2nd ed., Krieger, Malabar, FL, 1988.

\bibitem{GV} B. Gustafsson, A. Vasil'ev, \it Complex and Potential Analysis 
in Hele-Shaw cells. \rm Book, preliminary version, 2004.

\bibitem{HKZ} H. Hedenmalm, B. Korenblum, K. Zhu, \it Theory of Bergman 
spaces. \rm Graduate Texts in Mathematics, {\bf199}. Springer-Verlag, 
New York, 2000. 

\bibitem{HO} H. Hedenmalm, A. Olofsson, \it Hele-Shaw flow on weakly 
hyperbolic surfaces. \rm Indiana Univ. Math. J., to appear. 

\bibitem{HS} H. Hedenmalm, S. Shimorin, \it Hele-Shaw flow on hyperbolic 
surfaces. \rm J. Math. Pures Appl. (9) {\bf 81}  (2002),  no. 3, 187--222.

\bibitem{Kurt} K. Johansson, \it On fluctuations of eigenvalues of random 
Hermitian matrices. \rm Duke Math. J. {\bf 91} (1998), 151--204. 

\bibitem{KiS} D. Kinderlehrer, G. Stampacchia, \it An introduction to
variational inequalities and their applications, \rm Academic Press, New York,
1980.

\bibitem{KKMWWZ} I. K. Kostov, I. Krichever, M. Mineev-Weinstein, P. B.
Wiegmann, A. Zabrodin, \it $\tau$-function for analytic curves. \rm 
{\sl Random matrix models and their applications}, 285--299,
Math. Sci. Res. Inst. Publ., {\bf40}, Cambridge Univ. Press, Cambridge,
2001. 

\bibitem{KMZ} I. Krichever, A. Marshakov, A. Zabrodin, \it Integrable
structure of the Dirichlet boundary problem in multiply-connected domains. 
\rm Commun. Math. Phys., to appear. 

\bibitem{MWZ} A. Marshakov, P. Wiegmann, A. Zabrodin, \it Integrable
structure of the Dirichlet boundary problem in two dimensions. \rm Comm. Math.
Phys. {\bf 227} (2002), no. 1, 131--153. 

\bibitem{Meh} M. L. Mehta, \it Random matrices. \rm Second edition. Academic 
Press, Inc., Boston, MA, 1991. 


\bibitem{SaTo} E. B. Saff, V. Totik, {\it Logarithmic potentials with 
external fields.} Appendix B by Thomas Bloom. Grundlehren der Mathematischen 
Wissenschaften {\bf316}. Springer-Verlag, Berlin, 1997. 

\bibitem{Sakai} M. Sakai, \it Regularity of a boundary having a Schwarz
function, \rm Acta Math. \bf166 \rm (1991), 263--297. 

\bibitem{Wieg} 3] P. B. Wiegmann, \it Aharonov-Bohm effect in the quantum 
Hall regime and Laplacian growth problems.  \rm Statistical field theories 
(Como, 2001),  337--349, NATO Sci. Ser. II Math. Phys. Chem., {\bf73}, 
Kluwer Acad. Publ., Dordrecht, 2002.

\end{thebibliography}
\end{document}